\documentclass{article}%
\usepackage{amsmath}
\usepackage{amsfonts}
\usepackage{amssymb}
\usepackage{graphicx}%
\providecommand{\U}[1]{\protect\rule{.1in}{.1in}}
\begin{document}

\title{The Measurement of Statistical Evidence as the Basis for Statistical Reasoning}
\author{Michael Evans\\Dept. of Statistical\ Sciences\\University of Toronto}
\date{}
\maketitle

\noindent This paper serves as the basis for a talk at the OBayes Conference
at the U. of Warwick, June 2019.\bigskip

\noindent\textit{Abstract}: There are various approaches to the problem of how
one is supposed to conduct a statistical analysis. Different analyses can lead
to contradictory conclusions in some problems so this is not a satisfactory
state of affairs. It seems that all approaches make reference to the evidence
in the data concerning questions of interest as a justification for the
methodology employed. It is fair to say, however, that none of the most
commonly used methodologies is absolutely explicit about how statistical
evidence is to be characterized and measured. We will discuss the general
problem of statistical reasoning and the development of a theory for this that
is based on being precise about statistical evidence. This will be shown to
lead to the resolution of a number of problems.\medskip

\bigskip

\section{Introduction}

There are a variety of approaches to conducting statistical analyses and most
seem well-motivated from an intuitive point--of-view. If, in any particular
context, these all led to approximately the same results, then this wouldn't
be a problem. Unfortunately, it is possible that mutually contradictory
results can be produced and the Jeffreys-Lindley paradox is a well-known
example of this phenomenon. Given that these situations can arise, and with
very simple problems, it can lead one to wonder if there are any true
foundations for the subject. With the increasing importance of, what we will
call here, \textit{statistical reasoning} in almost all areas of science, this
ambiguity can be seen as an important problem to resolve.

It is interesting to note that virtually all approaches make reference to the
\textquotedblleft evidence in the data\textquotedblright\ or even to the
concept of \textquotedblleft statistical evidence\textquotedblright\ itself.
For example, p-values are considered as measuring the evidence against a
hypothesis although serious doubts have been raised about their suitability in
this regard. Similarly, likelihood ratios are considered as measuring the
evidence in the data supporting the truth of one value versus supporting the
truth of another and the Bayes factor is considered as a measure of
statistical evidence. There are also treatments that recognize the concept of
statistical evidence as a central aspect of statistical reasoning as with
Birnbaum (1964), Shafer (1976), Royall (1997), Thompson (2007),
Aitkin\ (2010), Morey, Romeijn, and Rouder (2016) and Vieland and Seok (2016).
There is also a long tradition of the consideration of evidence in the
philosophy of science, sometimes called confirmation theory, with Salmon
(1973) being a very accessible summary and Achinstein (2001) being a more
recent treatment.

While much of this literature has many relevant things to say about
statistical evidence, it is fair to say that there is no treatment that gives
an unambiguous definition of what statistical evidence is together with
developing a satisfactory theory of statistical reasoning based on this. Evans
(2015) contains a summary of a number of publications with co-authors that
have attempted to deal with this problem. This paper motivates and describes
that theory and provides a discussion of its advantages as well as its limitations.

As far as it benefits go, the theory provides a resolution of a number of
problems and this is illustrated via specific examples such as the
Jeffrey-Lindley paradox. It is also noted that the theory can be seen to
contribute to a resolution of a well-known divide in the statistical community
between frequentism and Bayesianism. In essence frequentism is concerned with
design: what are the consequences of the possible data sets that could be
obtained and how do we control for these effects? By contrast Bayesianism is
concerned with inference: given the observed data, what are the answers to
specific questions posed? So frequentism arises when the merits of a study are
being assessed and this surely has an impact on the acceptability of the
inferences produced via a Bayesian analysis.

While there are achievements of the approach taken here to statistical
reasoning, there are also limitations. Simply put, it is not the case that the
theory produces answers that are acceptable in all contexts where one might
deem statistical reasoning as appropriate. But such situations typically fall
outside the boundaries of what we might consider as acceptable statistical
problems and so require various compromises. For example, suppose one is
required to estimate $(\mu,\sigma^{2})$ based on the model $\{N(\mu,\sigma
^{2}):\mu\in R^{1},\sigma^{2}>0\}$ and a sample $x$ of size $n=1$ is
provided$.$ The estimate produced here is $(x,0)$ and one could argue that
this is absurd unless it is known categorically that there is no variability
in performances. Clearly, however, this can't be considered as a context where
the success of a theory is to be assessed. In such a situation some kind of
compromise is required, e.g., basing the estimate of $\sigma^{2}$ on \textit{a
prori} beliefs about this quantity. Such compromises can be seen as violating
some very basic principle of sound statistical reasoning, e.g., falsifiability
of the model and prior here. This points to the need to be careful when
identifying what the core problems are to which the theory is to be applied.
Otherwise, the possibility of obtaining an acceptable theory seems limited and
we would be left with the impression that the subject lacks a foundation. This
isn't an argument against compromises, but only that these are made when
necessary and that such compromises be clearly identified. In some ways this
is similar in spirit to the distinction made by Mosteller and Tukey (1977)
between exploratory and confirmatory analyses. The concern here is definitely
with the confirmatory aspect and moreover situations where design plays a role
through ensuring data is collected properly and an adequate amount is
collected according to various criteria subsequently discussed. Even further
restrictions are undoubtedly necessary to identify core statistical problems
and some discussion on this will be provided. It is worth noting here too the
papers Brown et al. (2002) and Ripamonti et al. (2017) that discuss
appropriate confidence and testing inferences for some binomial contexts. So
it is surely not correct to suggest that such problems are well-settled even
in the simplest of problems.\ When compromises are required it seems natural
to try and stay as close to a \textquotedblleft gold
standard\textquotedblright\ as possible, so it is necessary to identify what
that is.

\section{The Purpose of a Theory of Statistical Reasoning}

To start it is necessary to ask: what is the purpose of statistics as a
subject or what is its role in science? This is identified here as providing a
theory that gives answers to two questions based on observed data $\not x  $
and some object of interest $\Psi.$

\begin{quote}
\textbf{E: }What is a suitable value $\psi(x)$ for $\Psi$ together with an
assessment of the accuracy of this estimate?

\textbf{H}: Is there evidence that a specified value $\psi_{0}$ for $\Psi$\ is
true or false and how strong is this evidence?
\end{quote}

\noindent The requirement is placed on such a theory that the answers to
\textbf{E} and \textbf{H} be based on a clear characterization of the
statistical evidence in the data relevant to $\Psi.$

If the data $x$ was enough to answer these questions unambiguously that would
undoubtedly be best but it seems need more is needed. There are two aspects of
this that need discussion, referenced hereafter as the \textit{Ingredients}
and \textit{Statistical Inference}. These are to be taken as playing the same
roles in a statistical argument as the Premises and Rules of Inference play in
a logical argument or proof. For a logical argument may be \textit{valid}, as
whenever the rules of inference, such as \textit{modus ponens}, have been
applied correctly to the premises to derive a conclusion, but it is only a
\textit{sound} argument when it is valid and the premises are consistent and
true. An important aspect of a logical argument, see Kneale and Kneale (1962),
is that these two aspects of a logical argument are not confounded. So the
rules of inference are not changed based upon the premises and further that,
in applying the rules of inference, one is not concerned with the truth of the
premises. The validity of the argument and the truth of the premises are
separate issues. A similar position is adopted here.

\section{The Ingredients}

There are various criteria that the ingredients necessary for a theory of
statistical reasoning should satisfy. A partial list is as follows.

\begin{quote}
I$_{\text{1}}$ The minimum number of ingredients are required for a
characterization of statistical evidence so that \textbf{E} and \textbf{H} can
be answered.

I$_{\text{2}}$. The chosen ingredients are such that the bias in these choices
can be assessed and controlled.

I$_{\text{3}}$. Each element specified needs to be falsifiable via the data.

$\vdots$
\end{quote}

\noindent It seems clear that in the majority of statistical problems the
ingredients are chosen and are not specified by the application. As such the
ingredients are subjective and that would seem to violate a fundamental aspect
of science, namely, the goal of objectivity. The appropriateness of this is
recognized here but rather as an ideal that is not strictly achievable but
approached by dealing with the ingredients appropriately. I$_{\text{2}}$
suggests that it is possible to choose the ingredients in such a way that they
bias the answers to \textbf{E} and \textbf{H}. This is indeed the case, but as
long as the bias can be assessed and controlled, this need not be of great
concern. I$_{\text{3}}$ is yet another way of dealing with the inherent
subjectivity in a statistical analysis. For indeed it is reasonable to argue
that the chosen ingredients are never correct. In essence the role of the
ingredients is as devices that allow for a characterization of evidence so
that the rules of inference can be applied to answer \textbf{E} and
\textbf{H}. Their \textquotedblleft correctness\textquotedblright\ is
irrelevant unless the choices made are contradicted via the objective, when
collected correctly, data $x.$ The gold standard requires falsifiability of
all ingredients to help place statistical reasoning firmly within the realm of
science. Compromises to this need to be explicitly noted and only used when necessary.

So for a theory of statistical reasoning we need to state what the ingredients
are, how to check for their reasonableness and specify the rules of
statistical inference. The following ingredients, beyond the data $x,$ are
necessary for the developments here.

\begin{quote}
\textbf{Ingredients}

Model:$\{f_{\theta}:\theta\in\Theta\}$ a collection of conditional probability
distributions for $x\in\mathcal{X}$ given $\theta$ such that the object of
interest $\psi=\Psi(\theta)$ is specified by the true distribution that gave
rise to $x.$

Prior: $\pi$ a probability distribution on $\Theta.$

Delta: $\delta$ the difference that matters so that dist$(\psi_{1},\psi
_{2})\leq\delta$, for some distance measure dist, means that $\psi_{1}$ and
$\psi_{2}$ are for practical purposes indistinguishable.
\end{quote}

\noindent So the model and prior specify a joint probability distribution for
$\omega=(\theta,x)\sim\pi(\theta)f_{\theta}(x).$ As such all uncertainties are
handled within the context of probability which is interpreted here as
measuring belief, no matter how the probabilities are assigned.

The role of $\delta$ will be subsequently discussed but it raises an
interesting and relevant point concerning the role of infinities in
statistical modelling. The position taken here is that whenever infinities
appear their role is as approximations as expressed via limits and they do not
represent reality. For example, data arises via a measurement process and as
such all data is measured to a finite accuracy.\ So data is discrete and
moreover sample spaces are bounded as measurements cannot be arbitrarily
large/small positively or negatively. So whenever continuous probability
distributions are used these are considered as approximations to essentially
finite distributions. There are a variety of opinions about the use of
infinity among physicists, mathematicians, statisticians, etc., but little is
lost by taking this view of things here and there is a substantial gain for
the theory through the avoidance of anomalies.

The question now is how to choose the model and the prior? Unfortunately, we
are somewhat silent about general principles for choosing a model, but when it
comes to the prior for us this is by an elicitation algorithm which explains
why the prior in question has been chosen. An inability to come up with a
suitable elicitation suggests a lack of sufficient understanding of the
scientific problem or an inappropriate choice of model where the real world
meaning of $\theta$ is somewhat unclear. The existence of a suitable
elicitation algorithm could be viewed as a necessity to place a context within
the gold standard but, given the way models are currently chosen, we do not
adopt that position. Still whatever approach is taken to choosing $\pi,$ it is
subject to I$_{\text{2}}$ and I$_{\text{3}}$. As will be seen, the
implementation of I$_{\text{2}}$ requires the characterization of statistical
evidence and so discussion of this is delayed until Section 6 and
I$_{\text{3}}$ is addressed next.

\section{Checking the Ingredients}

If the observed data $x$ is surprising (in the \textquotedblleft
tails\textquotedblright\ of) for each distribution in $\{f_{\theta}:\theta
\in\Theta\},$ then this suggests a problem with the model, and otherwise the
model is at least acceptable. There are a number of approaches available for
checking the model and this isn't discussed further here although the model is
undoubtedly the most important ingredient chosen.

To be a bit more formal note that, when $T$ is \ a minimal sufficient
statistic for the model, then the joint factors as
\[
\pi(\theta)f_{\theta}(x)=\pi(\theta)f_{\theta,T}(T(x))f(x\,|\,T(x))=\pi
(\theta\,|\,T(x)m_{T}(T(x))f(x\,|\,T(x))
\]
where $f(\cdot\,|\,T(x))$ is a probability distribution, independent of
$\theta,$ available for model checking, $m_{T}$ is the prior (predictive)
distribution of $T,$ available for checking the prior, while $\pi
(\cdot\,|\,T(x))$ is the posterior of $\theta$\ and provides probabilities for
the inference step. Note that Box (1980) proposed using the prior (predictive)
distribution of $x$ for jointly checking the model and prior but, to ascertain
where a problem exists when it does, it is more appropriate to split this into
checking the model first and then, if the model passes its checks, check the prior.

A prior fails when the true value lies in its tails. Evans and Moshonov (2006)
proposed using the tail probability
\begin{equation}
M_{T}(m_{T}(t)\leq m_{T}(T(x))) \label{eq1}%
\end{equation}
for this purpose. This was generalized to include conditioning on ancillary
statistics, to remove variation irrelevant to checking the prior, and also to
further factorizations of $m_{T}(T(x))$ that allow for checking of individual
components of the prior, to try and isolate a problem when one exists. In
Evans and Jang (2011a) it is established that, under conditions, $M_{T}%
(m_{T}(t)\leq m_{T}(T(x)))\rightarrow\Pi(\pi(\theta)\leq\pi(\theta_{true}))$
as the amount of data increases and so the tail probability is a consistent
check on the prior. Further refinements can be proposed to deal with
invariance and Nott et al. (2018) generalizes (\ref{eq1}) to provide a fully
invariant check that connects nicely with the measure of evidence used for
inference. It is shown in Al Labadi and Evans (2017) that, when prior-data
conflict exists, then inferences can be very sensitive to perturbations of the prior.

What does one do when a prior fails? Evans and Jang (2011b) provides an
approach to defining what is meant by weakly informative priors based on an
idea in Gelman et al. (2008). In essence a definition is provided for what it
means for one prior to be weakly informative with respect to another and
further quantifies this in terms of fewer prior-data conflicts expected. This
does not involve the observed data. One can conceive of a base elicited prior
and a hierarchy of successively more weakly informative priors. So a prior for
which conflict is detected, can be replaced by one more weakly informative and
one can progress up the hierarchy until conflict is avoided. The point here is
that final prior is not strictly speaking dependent on the data.

There are many approaches to prior-data conflict discussed in the literature
and this seems to be, like model checking, the kind of problem for which there
isn't a definitive methodology. The significance of the work cited, however,
is that it presents a package for dealing with the reasonableness of the prior
as part of a theory of statistical reasoning.

\section{The Rules of Statistical Inference}

There are three rules of inference that are used to determine inferences.
These are stated for a probability model $(\Omega,\mathcal{F},P).$ Suppose
interest is in whether or not the event $A\in\mathcal{F}$ is true after
observing $C\in\mathcal{F}$ and it is supposed that both $P(A)>0$ and $P(C)>0$
with the null case handled via taking appropriate limits$.$

\begin{quote}
R$_{\text{1}}.$ \textit{Principle of conditional probability}: beliefs about
$A,$ as expressed initially by $P(A),$ are replaced by $P(A\,|\,C).$

R$_{\text{2}}$. \textit{Principle of evidence}: the observation of $C$ is
evidence in favor of $A$ when $P(A\,|\,C)>P(A),$ is evidence against $A$ when
$P(A\,|\,C)<P(A)$ and is evidence neither for nor against $A$ when
$P(A\,|\,C)=P(A).$

R$_{\text{3}}$. The evidence is measured quantitatively by the relative belief
ratio%
\[
RB(A\,|\,C)=\frac{P(A\,|\,C)}{P(A)}.
\]

\end{quote}

\noindent While R$_{\text{1}}$ doesn't seem controversial, its strict
implementation in Bayesian contexts demands proper priors and priors that do
not depend on the data. R$_{\text{2}}$ also seems quite natural and, as will
be seen, really represents the central core of our approach to statistical
reasoning. R$_{\text{3}}$ is perhaps not as obvious, but it is clear that
$RB(A\,|\,C)>1(<,=)$ indicates that evidence in favor of (against, neither)
$A$ has been obtained. In fact the relative belief ratio only plays a role
when it is necessary to order alternatives.

There are other \textit{valid} measures of evidence in the sense that they
have a cut-off that determines evidence in favor or against, as specified by
R$_{\text{2}}$, and can be used to order alternatives, see the discussion in
Evans (2015). The relative belief ratio, however, has some advantages, such as
invariance under reparameterizations when continuous models are being used and
see Section 7. As will be seen, any increasing function of $RB$ leads to
identical inferences. Other possible valid measures include the difference
$P(A\,|\,C)-P(A),$ with cut-off 0, or the Bayes factor
\[
BF(A\,|\,C)=\frac{P(A\,|\,C)}{P(A^{c}\,|\,C)}/\frac{P(A)}{P(A^{c})}%
=\frac{RB(A\,|\,C)}{RB(A^{c}\,|\,C)},
\]
with cut-off 1. As indicated, the Bayes factor can be defined in terms of the
relative belief ratio\ but not conversely. Since $RB(A\,|\,C)>1$ iff
$RB(A^{c}\,|\,C)<1,$ the Bayes factor isn't really a comparison of the
evidence for $A$ with the evidence for its negation. Furthermore, in the
continuous case, when $RB$ and $BF$\ are defined as limits as below, they are
equal. So $RB$ is preferred to $BF.$

Now consider the Bayesian context. When interest is in $\psi=\Psi(\theta),$
then generally the relative belief ratio at a value $\psi$ equals
\[
RB_{\Psi}(\psi\,|\,x)=\lim_{\epsilon\rightarrow0}RB(N_{\epsilon}%
(\psi)\,|\,x)=\frac{\pi_{\Psi}(\psi\,|\,x)}{\pi_{\Psi}(\psi)}%
\]
where neighborhoods $N_{\epsilon}(\psi)$ of $\psi$ satisfy $N_{\epsilon}%
(\psi)\overset{\text{nicely}}{\rightarrow}\{\psi\}$ as $\epsilon\rightarrow0$
and the equality on the right holds whenever the prior density $\pi_{\Psi}$ of
$\Psi$ is positive and continuous at $\psi.$ This definition is motivated by
the treatment of contexts where infinities arise.

\subsection{Problem \textbf{E}}

Suppose that the range of possible values for $\psi=\Psi(\theta)$ is also
denoted by $\Psi.$ Then R$_{\text{3}}$ determines the relative belief estimate
as
\[
\psi(x)=\arg\sup_{\psi\in\Psi}RB_{\Psi}(\psi\,|\,x)
\]
as this maximizes the evidence in favor and note $\sup_{\psi\in\Psi}RB_{\Psi
}(\psi\,|\,x)\geq1\ $always and generally the inequality is strict. To measure
the accuracy of $\psi(x)$ there are a number of possibilities but the
\textit{plausible region}
\[
Pl_{\Psi}(x)=\{\psi:RB_{\Psi}(\psi\,|\,x)>1\},
\]
the set of $\psi$ values where there is evidence in favor of the value being
true, is surely central to this. If the \textquotedblleft
size\textquotedblright, such as volume or prior content, of $Pl_{\Psi}(x)$ is
small and its posterior content $\Pi_{\Psi}(Pl_{\Psi}(x)\,|\,x)$ is large,
then this suggests an accurate estimate has been obtained.

There are several notable aspects of this. First, in the continuous case, the
methodology is invariant. So if instead interest is in $\lambda=\Lambda
(\Psi(\theta)),$ where $\Lambda$ is 1-1 and smooth, then $\lambda
(x)=\Lambda(\psi(x))$ and $Pl_{\Lambda}(x)=\Lambda(Pl_{\Psi}(x)).$ Second,
while $\psi(x)$ can also be thought of as the MLE from an integrated
likelihood, that approach does not lead to $Pl_{\Psi}(x)$ because likelihoods
do not define evidence for specific values. Third, and perhaps most
significant, the set $Pl_{\Psi}(x)$ is completely independent of how evidence
is measured quantitatively. In other words if any valid measure of evidence is
used, then the same set $Pl_{\Psi}(x)$ is obtained. This has the happy
consequence that, however we choose to estimate $\Psi$ via an estimator that
respects the principle of evidence, then effectively the same quantification
of error is obtained. This points to the possibility of using some kind of
smoothing operation on $\psi(x)$ to produce values that lie in $Pl_{\Psi}(x)$
when this is considered necessary$.$ It is also possible to use, as part of
measuring the accuracy of $\psi(x),$ a $\gamma$-relative belief credible
region for $\Psi,$ namely,%
\[
C_{\Psi,\gamma}(x)=\{\psi:RB_{\Psi}(\psi\,|\,x)\geq c_{\Psi,\gamma}(x)\}
\]
where $c_{\Psi,\gamma}(x)=\inf_{c}\{c:\Pi_{\Psi}(\{\psi:RB_{\Psi}%
(\psi\,|\,x)>c\}\,|\,x)<\gamma\}.$ To use this region, however, it is
necessary that $\gamma\leq\Pi_{\Psi}(Pl_{\Psi}(x)\,|\,x)$ otherwise
$C_{\Psi,\gamma}(x)$ contains values for which there is evidence against and
so would not accurately represent the evidence in the data.

The following example illustrates that probabilities do not measure evidence,
although this is a common misconception.\smallskip

\noindent\textbf{Example 1.} \textit{Prosecutor's Fallacy}

Assume a uniform probability distribution on a population of size $N$ of which
some member has committed a crime. DNA evidence has been left at the crime
scene and suppose this trait is shared by $m\ll N$ of the population. A
prosecutor is criticized because they conclude that, because the trait is rare
and a particular member possesses the trait, then they are guilty. In fact
they misinterpret $P($\textquotedblleft has trait\textquotedblright%
$\,|\,\text{\textquotedblleft}$guilty$\text{\textquotedblright})=1$ as the
probability of guilt which is $P($\textquotedblleft guilty\textquotedblright%
$\,|\,$"has trait\textquotedblright$)=1/m$ which is small if $m$ is large. But
this probability does not reflect the evidence of guilt. For, if you have the
trait, then clearly this is evidence in favor of guilt. Note that
\begin{align*}
RB(\text{\textquotedblleft guilty\textquotedblright}\,|\,\text{"has
trait\textquotedblright}))  &  =\frac{P(\text{\textquotedblleft
guilty\textquotedblright}\,|\,\text{\textquotedblleft has
trait\textquotedblright})}{P(\text{\textquotedblleft guilty\textquotedblright%
)}}=\frac{1/m}{1/N}=\frac{N}{m}>1\\
RB(\text{\textquotedblleft not guilty\textquotedblright}\,|\,\text{"has
trait\textquotedblright}))  &  =\frac{P(\text{\textquotedblleft not
guilty\textquotedblright}\,|\,\text{\textquotedblleft has
trait\textquotedblright})}{P(\text{\textquotedblleft not
guilty\textquotedblright)}}\\
&  =\frac{(m-1)/m}{(N-1)/N}=\frac{N}{N-1}\frac{m}{m-1}<1
\end{align*}
and $Pl($\textquotedblleft has trait\textquotedblright$)=\{$\textquotedblleft
guilty\textquotedblright$\}$ with posterior content $1/m.$ So there
\textit{is} evidence of guilt but it is weak whenever $m$ is large. Note too
that the MAP estimate is \textquotedblleft not guilty\textquotedblright$\,$.
So the prosecutor is correct that there is evidence of guilt and because, in
general, $RB(A\,|\,C)=RB(C\,|\,A)$ this conclusion would have been supported
if they had used $P($"has trait\textquotedblright$\,|\,\text{\textquotedblleft%
}$guilty$\text{\textquotedblright})$ correctly as part of a relative belief
ratio but of course the strength measurement would have been
incorrect.\smallskip

\noindent A general weakness of MAP/HPD inferences is their lack if invariance
under reparameterizations. Example 1 shows, however, that these inferences are
generally inappropriate because they do not express evidence properly. Example
1 also demonstrates a distinction between decisions and inferences. Clearly
when $m$ is large there should not be a conviction on the basis of weak
evidence. But suppose that \textquotedblleft guilty\textquotedblright%
\ corresponds to being a carrier of a highly infectious deadly disease and
\textquotedblleft has trait\textquotedblright\ corresponds to some positive
(but not definitive) test for this, then the same numbers should undoubtedly
lead to a quarantine. In essence a theory of statistical reasoning should tell
us what the evidence says and decisions are made, partly on this basis, but
employing many other criteria as well.

\subsection{Problem \textbf{H}}

It is immediate that $RB_{\Psi}(\psi_{0}\,|\,x)$ is the evidence concerning
$H_{0}:\Psi(\theta)=\psi_{0}.$ The evidential ordering implies that the
smaller $RB_{\Psi}(\psi_{0}\,|\,x)$ is than 1, the stronger the evidence is
against $H_{0}$ and the bigger it is than 1, the stronger the evidence is in
favor $H_{0}.$ But how is one to measure this strength? In Baskurt and Evans
(2013) it is proposed to measure the \textit{strength of the evidence} via%
\begin{equation}
\Pi_{\Psi}\left(  \left.  RB_{\Psi}(\psi\,|\,x)\leq RB_{\Psi}(\psi
_{0}\,|\,x)\,\right\vert \,x\right)  \label{eq2}%
\end{equation}
which is the posterior probability that the true value of $\psi$ has evidence
no greater than that obtained for the hypothesized value $\psi_{0}.$ When
$RB_{\Psi}(\psi_{0}\,|\,x)<1$ and (\ref{eq2}) is small, then there is strong
evidence against $H_{0}$ since there is a large posterior probability that the
true value of $\psi$ has a larger relative belief ratio. Similarly, if
$RB_{\Psi}(\psi_{0}\,|\,x)>1$ and (\ref{eq2}) is large, then there is strong
evidence that the true value of $\psi$ is given by $\psi_{0}$ since there is a
large posterior probability that the true value is in $\{\psi:RB_{\Psi}%
(\psi\,|\,x)\leq RB_{\Psi}(\psi_{0}\,|\,x)\}$ and $\psi_{0}$ maximizes the
evidence in this set. The strength measurement here results from comparing the
evidence for $\psi_{0}$ with the evidence for each of the other possible
$\psi$ values. This is appropriate provided the values are all referencing
similar objects. Note that, using Markov's inequality,
\[
\Pi_{\Psi}\left(  \{\psi_{0}\}\,|\,x\right)  \leq\Pi_{\Psi}\left(  \left.
RB_{\Psi}(\psi\,|\,x)\leq RB_{\Psi}(\psi_{0}\,|\,x)\,\right\vert \,x\right)
\leq RB_{\Psi}(\psi_{0}\,|\,x)
\]
and so a very small value of $RB_{\Psi}(\psi_{0}\,|\,x)$ is immediately strong
evidence against $H_{0}.$ Also, when there is evidence in favor, then
$\psi_{0}\in Pl_{\Psi}(x)$ and so the size and posterior content of this set
also provides an indication of the strength of the evidence.

When $H_{0}$ is false, then $RB_{\Psi}(\psi_{0}\,|\,x)$ converges to 0 as does
(\ref{eq2}). When $H_{0}$ is true, then $RB_{\Psi}(\psi_{0}\,|\,x)$ converges
to its largest possible value (greater than 1 and often $\infty$) and, in the
discrete case (\ref{eq2}) converges to 1. In the continuous case, however,
when $H_{0}$ is true, then (\ref{eq2}) typically converges to a $U(0,1)$
random variable. This is easily resolved by requiring that a deviation
$\delta>0$ be specified such that if dist$(\psi_{1},\psi_{2})<\delta,$ where
dist is some measure of distance determined by the application, then this
difference is to be regarded as immaterial. This leads to redefining $H_{0}$
as $H_{0}=\{\psi:$ dist$(\psi,\psi_{0})<\delta\}$ and typically a natural
discretization of $\Psi$ exists with\ $H_{0}$ as one of its elements. With
this modification (\ref{eq2}) converges to 1 as the amount of data increases
when $H_{0}$ is true. Some discussion on determining a relevant $\delta$ can
be found in Al-Labadi, Baskurt and Evans (2017) and Evans, Guttman and Li
(2017). Typically the incorporation of such a $\delta$ makes computations
easier as then there is no need to estimate densities when these are not
available in closed form.

The following simple example demonstrates the need to calibrate the measure of
evidence.\smallskip

\noindent\textbf{Example 2. }\textit{Location normal.}

Suppose $x=(x_{1},\ldots,x_{n})$ is i.i.d.\ $N(\mu,\sigma_{0}^{2})$ with
$\sigma_{0}^{2}$ known and $\pi$ is a $N(\mu_{\ast},\tau_{\ast}^{2})$ prior
and the hypothesis is $H_{0}:\mu=\mu_{0}.$ So%
\begin{align*}
&  RB(\mu_{0}\,|\,x)\\
&  =\left(  1+\frac{n\tau_{\ast}^{2}}{\sigma_{0}^{2}}\right)  ^{1/2}%
\exp\left\{
\begin{array}
[c]{c}%
-\frac{1}{2}\left(  1+\frac{\sigma_{0}^{2}}{n\tau_{\ast}^{2}}\right)
^{-1}\left(  \frac{\sqrt{n}(\bar{x}-\mu_{0})}{\sigma_{0}}+\frac{\sigma_{0}%
(\mu_{\ast}-\mu_{0})}{\sqrt{n}\tau_{\ast}^{2}}\right)  ^{2}\\
+\frac{\left(  \mu_{0}-\mu_{\ast}\right)  ^{2}}{2\tau_{0}^{2}}%
\end{array}
\right\}  ,
\end{align*}
which, in this case is the same as the Bayes factor for $\mu_{0}$\ obtained
via Jeffreys' mixture approach. From this it is easy to see that $RB(\mu
_{0}\,|\,x)\rightarrow\infty$ as $\tau_{\ast}^{2}\rightarrow\infty$ or, when
$\sqrt{n}|\bar{x}-\mu_{0}|/\sigma_{0}$ is fixed, as $n\rightarrow\infty.$ If
this is calibrated using the strength then, under these circumstances%
\begin{equation}
\Pi\left(  \left.  RB(\mu\,|\,x)\leq RB(\mu_{0}\,|\,x)\,\right\vert
\,x\right)  \rightarrow2(1-\Phi(\sqrt{n}|\bar{x}-\mu_{0}|/\sigma_{0}))
\label{eq3}%
\end{equation}
the classical p-value for assessing $H_{0}.$ The Jeffreys-Lindley paradox is
the apparent divergence between $RB(\mu_{0}\,|\,x)$ and the p-value as
measures of evidence concerning $H_{0}.$ For example, see Robert (2014), Villa
and Walker (2017) for some recent discussion of the paradox$.$ Using the
relative belief framework, however, it is clear that while $RB(\mu_{0}\,|\,x)$
may be a large value, and so apparently strong evidence in favor, (\ref{eq3})
suggests it is weak evidence in favor. This doesn't fully explain why this
arises, however, as clearly when $\sqrt{n}|\bar{x}-\mu_{0}|/\sigma_{0}$ is
large, then we would expect there to be evidence against. To understand why
this discrepancy arises, it is necessary to consider bias as discussed in the
next section.

Another aspect of (\ref{eq3}) is of interest because it raises the obvious
question: is the classical p-value a valid measure of evidence? The answer is
no according to our definition, but the difference of two tail probabilities
\begin{align*}
&  2\left(  1-\Phi\left(  \frac{\sqrt{n}|\bar{x}-\mu_{0}|}{\sigma_{0}}\right)
\right)  -\\
&  2\left(  1-\Phi\left(  \left[  \log\left(  1+\frac{n\tau_{\ast}^{2}}%
{\sigma_{0}^{2}}\right)  )+\left(  1+\frac{n\tau_{\ast}^{2}}{\sigma_{0}^{2}%
}\right)  ^{-1}\frac{\left(  \bar{x}-\mu_{\ast}\right)  ^{2}}{\tau_{0}^{2}%
}\right]  ^{1/2}\right)  \right)  ,
\end{align*}
with cut-off 0, is a valid measure of evidence. Note that the right-most tail
probability converges almost surely to 0 as $n\rightarrow\infty$ or $\tau
_{0}^{2}\rightarrow\infty$ and this implies that, if the classical p-value is
to be used to measure evidence, the significance level has to go to 0 as $n$
increases. Note that if the $N(\mu_{\ast},\tau_{\ast}^{2})$ prior is replaced
by a Uniform$(-m,m)$ prior, then the same conclusion is reached as
$n\rightarrow\infty$ or as $m\rightarrow\infty.$ These conclusions are similar
to those found in Berger and Selke (1987) and Berger and Delampady (1987).

\section{Measuring Bias}

One of the more serious concerns with Bayesian methodology is with the issue
of bias. Can an analyst choose the ingredients in such a way that the results
are a foregone conclusion with high prior probability? The answer is yes and
the situation described in Example 1 is a good example of this. Still it is
necessary to be precise about the meaning of bias and the principle of
evidence provides this. Even though we will use the relative belief ratio in
the discussion here, it is important to note that the bias measures discussed
are the same no matter what valid measure of evidence is used. The approach to
measuring bias described here for problem \textbf{H} appeared in Baskurt and
Evans (2013) and this has since been extended to problem \textbf{E }in Evans
and Guo (2019). Bias should be measured a priori as it can be controlled by
design although a post hoc measurement is also possible. In essence the bias
numbers give a measure of the merit of a statistical study. Studies with high
bias cannot be considered as being reliable irrespective of the correctness of
the ingredients. The measurement of bias also establishes important links with frequentism.

\subsection{Bias for Problem \textbf{H}}

Consider the hypothesis $H_{0}:\Psi(\theta)=\psi_{0}$ and let $M(\cdot
\,|\,\psi)$ denote the prior predictive distribution of the data given that
$\Psi(\theta)=\psi.$ In general it is possible for there to be a high prior
probability that evidence against $H_{0}$ will be obtained even when it is
true. In such a circumstance, if evidence against $H_{0}$ is obtained, then
the result seems highly suspect.\ Similarly, if there is a high prior
probability of obtaining evidence in favor of $H_{0}$ even when it
meaningfully false, then actually obtaining such evidence based on observed
data is hardly compelling. So measuring bias against and bias in favor of
$H_{0}$ are important aspects of this problem.

Bias against $H_{0}$ means that the ingredients and the data $x,$ of a
prescribed size, are such that with high prior probability evidence will not
be obtained in favor of $H_{0}$ even when it is true. Bias against $H_{0}$ is
thus measured by
\begin{equation}
M(RB_{\Psi}(\psi_{0}\,|\,x)\leq1\,|\,\psi_{0}). \label{biasagainst}%
\end{equation}
If (\ref{biasagainst}) is large, then obtaining evidence against $H_{0}$ seems
like a foregone conclusion. For bias in favor, consider a distance measure
dist :$\Psi^{2}\rightarrow\lbrack0,\infty).$ Bias in favor of $H_{0}$ is
measured by
\begin{equation}
\sup_{\psi_{\ast}\in\{\psi:\text{ dist}(\psi,\psi_{0})\geq\delta\}}M(RB_{\Psi
}(\psi_{0}\,|\,x)\geq1\,|\,\psi_{\ast}) \label{biasfor}%
\end{equation}
and note the necessity of including $\delta$ so that the measure is based only
on those values of $\psi$ that are meaningfully different from $\psi_{0}.$
Typically $M(RB_{\Psi}(\psi_{0}\,|\,D)\geq1\,|\,\psi_{\ast})$ increases as
dist$(\psi_{\ast},\psi_{0})$ increases so the supremum can then be taken
instead over $\{\psi:$ dist$(\psi,\psi_{0})=\delta\}.$ When (\ref{biasfor}) is
large, then it suggests that finding evidence in favor of $H_{0}$ is
meaningless. The choice of the prior can be used somewhat to control bias but
typically a prior that makes one bias lower just results in making the other
bias higher. It is established in Evans (2015) that, under quite general
circumstances, both biases converge to 0 as the amount of data increases. So
bias can be controlled by design a priori.

The following example illustrates some general characteristics of measuring
and controlling biases.\smallskip

\noindent\textbf{Example 3. }\textit{Location normal (continued).}

The bias against and bias in favor can be computed in closed form in this
case, see Evans and Guo (2019). Table 1 gives some values for the bias against
when testing the hypothesis $H_{0}:\mu=0$ when $\sigma_{0}^{2}=1$ for two
priors$.$ The results here illustrate something that holds generally in this
case. Provided the prior variance $\tau_{\ast}^{2}>\sigma_{0}^{2}/n,$ the
maximum bias against is achieved when the prior mean $\mu_{\ast}$ equals the
hypothesized mean $\mu_{0}.$ This turns out to be very convenient\ as
controlling the bias against for this case controls the bias against
everywhere. This seems paradoxical, as the maximum amount of belief is being
placed at $\mu_{0}$ when $\mu_{\ast}=\mu_{0}.$ It is clear, however, that the
more prior probability that is assigned to a value the harder it is for the
probability to increase. This is another example of a situation where evidence
works somewhat contrary to our intuition based on how belief works. Another
conclusion from Table 2 is that bias against is not a problem.\ Provided the
prior isn't chosen too concentrated, this is generally the case, at least in
our experience.
%TCIMACRO{\TeXButton{B}{\begin{table}[tbp] \centering}}%
%BeginExpansion
\begin{table}[tbp] \centering
%EndExpansion%
\begin{tabular}
[c]{|c|c|c|}\hline
$n$ & $(\mu_{\ast},\tau_{\ast}^{2})=(1,1)$ & $(\mu_{\ast},\tau_{\ast}%
^{2})=(0,1)$\\\hline
$5$ & $0.095$ & $0.143$\\
$10$ & $0.065$ & $0.104$\\
$20$ & $0.044$ & $0.074$\\
$50$ & $0.026$ & $0.045$\\
$100$ & $0.018$ & $0.031$\\\hline
\end{tabular}
\caption{Bias against for the hypothesis
$H_0=\{0\}$ with a $N(\mu_*,\tau_*^2)$ prior for different sample sizes $n$ with $\sigma_0=1$.}\label{loctab1}%
%TCIMACRO{\TeXButton{E}{\end{table}}}%
%BeginExpansion
\end{table}%
%EndExpansion

\begin{figure}
[ptb]
\begin{center}
\includegraphics[
height=2.7069in,
width=2.7069in
]%
{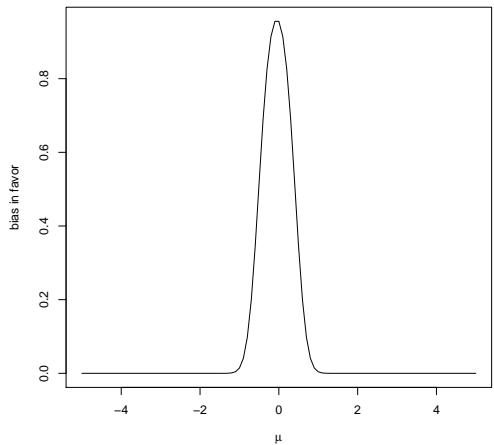}%
\caption{Plot of $M(RB(0\,|\,X)\geq1\,|\,\mu)$ when $n=20,\mu_{\ast}%
=1,\tau_{\ast}=1,\sigma_{0}=1.$}%
\label{fig1}%
\end{center}
\end{figure}

Figure \ref{fig1} is a plot of $M(RB(0\,|\,x)\geq1\,|\,\mu)$ as a function of
$\mu$ for this problem when $n=20,\mu_{\ast}=1,\tau_{\ast}^{2}=1.$ This
demonstrates how the the bias in favor drops off as the true value moves away
from $\mu_{0}.$ Table 2 contains values of the bias in favor for several
priors with $\delta=0.5$ and dist taken to be Euclidean distance. It is clear
that bias in favor is a much more serious issue. To a certain extent this can
be reduced by choosing $\delta$ larger, but $\delta$ is determined by the
application and really only sample size is available to control bias.

As discussed in Evans and Guo (2019), when $\tau_{\ast}^{2}\rightarrow\infty$
the bias against goes to 0 and the bias in favor goes to 1 in this problem.
This explains the phenomenon associated with the Jeffreys-Lindley paradox, as
taking a very diffuse prior induces extreme bias in favor. So this argues
against using arbitrarily diffuse priors.\ Rather one should elicit the value
for $(\mu_{\ast},\tau_{\ast}^{2})$, prescribe the value of $\delta$ and choose
$n$ to make the biases acceptable.%
%TCIMACRO{\FRAME{ftbpFU}{2.7069in}{2.7069in}{0pt}{\Qcb{Plot of
%$M(RB(0\,|\,X)\geq1\,|\,\mu)$ when $n=20,\mu_{\ast}=1,\tau_{\ast}=1,\sigma
%_{0}=1.$}}{\Qlb{fig1}}{fig2.eps}{\special{ language "Scientific Word";
%type "GRAPHIC";  maintain-aspect-ratio TRUE;  display "USEDEF";
%valid_file "F";  width 2.7069in;  height 2.7069in;  depth 0pt;
%original-width 6.9998in;  original-height 6.9998in;  cropleft "0";
%filename 'fig2.eps';file-properties "XNPEU";}}}%
%BeginExpansion
%EndExpansion
%

%TCIMACRO{\TeXButton{B}{\begin{table}[tbp] \centering}}%
%BeginExpansion
\begin{table}[tbp] \centering
%EndExpansion%
\begin{tabular}
[c]{|c|c|c|}\hline
$n$ & $(\mu_{\ast},\tau_{\ast}^{2})=(1,1)$ & $(\mu_{\ast},\tau_{\ast}%
^{2})=(0,1)$\\\hline
$5$ & $0.871$ & $0.631$\\
$10$ & $0.747$ & $0.516$\\
$20$ & $0.519$ & $0.327$\\
$50$ & $0.125$ & $0.062$\\
$100$ & $0.006$ & $0.002$\\\hline
\end{tabular}
\caption{Bias in favor of the hypothesis
$H_0=\{0\}$ with a $N(\mu_*,\tau_*^2)$ prior  for different sample sizes $n$ with $\sigma_0=1$ and $\delta=0.5$.}\label{loctab2}%
%TCIMACRO{\TeXButton{E}{\end{table}}}%
%BeginExpansion
\end{table}%
%EndExpansion

\subsection{Bias for Problem \textbf{E}}

The accuracy of a valid estimate of $\psi$ is measured by the size of the
plausible region $Pl_{\Psi}(x)=\{\psi:RB_{\Psi}(\psi\,|\,x)>1\}$. As such, if
the plausible region is reported as containing the true value and it does not,
then the evidence is misleading. It is argued therefore, that biases in
estimation problems be measured by prior coverage probabilities associated
with $Pl_{\Psi}(x)$ and the implausible region $Im_{\Psi}(x)=\{\psi:RB_{\Psi
}(\psi\,|\,x)<1\},$ the set of values for which there is evidence against. As
will be shown, there is a close connection between measuring bias for problems
\textbf{H} and \textbf{E}.

The prior probability that the plausible region does not cover the true value
measures bias against when estimating $\psi.$ This equals
\begin{equation}
E_{\Pi_{\Psi}}\left(  M(\psi\notin Pl_{\Psi}(X)\,|\,\psi)\right)
=E_{\Pi_{\Psi}}(M(RB_{\Psi}(\psi\,|\,X)\leq1\,|\,\psi)), \label{biasagest1}%
\end{equation}
which is also the average bias against over all hypothesis testing problems
$H_{0}:\Psi(\theta)=\psi.$ Note that $1-E_{\Pi_{\Psi}}\left(  M(\psi\notin
Pl_{\Psi}(X)\,|\,\psi)\right)  =E_{\Pi_{\Psi}}\left(  M(\psi\in Pl_{\Psi
}(X)\,|\,\psi)\right)  $\newline$=E_{M}\left(  \Pi_{\Psi}(Pl_{\Psi
}(X)\,|\,X)\right)  $ which is the prior coverage probability of $Pl_{\Psi}.$
So, if the bias against for estimation is small, the coverage probability is
high$.$ Note that
\begin{equation}
\sup_{\psi}M(\psi\notin Pl_{\Psi}(X)\,|\,\psi)=\sup_{\psi}M(RB_{\Psi}%
(\psi\,|\,X)\leq1\,|\,\psi), \label{biasagest2}%
\end{equation}
is an upper bound on (\ref{biasagest1}). Therefore, controlling
(\ref{biasagest2}) controls the bias against in estimation and all hypothesis
assessment problems involving $\psi$. Also $1-\sup_{\psi}M(\psi\notin
Pl_{\Psi}(X)\,|\,\psi)=\inf_{\psi}M(\psi\in Pl_{\Psi}(X)\,|\,\psi)\leq
E_{M}\left(  \Pi_{\Psi}(Pl_{\Psi}(X)\,|\,X)\right)  $ so using
(\ref{biasagest2}) implies lower bounds for the coverage probability of the
plausible region and for the expected posterior content of the plausible
region. In general, both (\ref{biasagest1}) and (\ref{biasagest2}) converge to
0 with increasing amounts of data. So it is possible to control for bias
against in estimation problems by design.

The prior coverage probability that follows from (\ref{biasagest1}) can be
considered as a Bayesian confidence as it is dependent on the full prior. The
lower bound on the coverage probability that arises from (\ref{biasagest2})
only depends on the conditional priors on the nuisance parameters and such a
situation could be considered as somewhat like the use of priors in random
effects models. When $\Psi(\theta)=\psi$ then the coverage probability is a
\textquotedblleft pure frequentist\textquotedblright\ coverage probability.
\smallskip\newpage

\noindent\textbf{Example 4. }\textit{Location normal (continued).}

Table \ref{loctab3} gives some values of the bias against measure for
estimation with different priors and sample sizes. Subtracting the
probabilities in Table \ref{loctab3} from 1 gives the prior probability that
the plausible region covers the true value and the expected posterior content
of the plausible region. So when $n=20,\tau_{\ast}=1,$ the prior probability
of the plausible region containing the true value is $1-0.051=0.949$ so
$Pl(x)$ is a $0.949$ Bayesian confidence interval for $\mu.$%

%TCIMACRO{\TeXButton{B}{\begin{table}[tbp] \centering}}%
%BeginExpansion
\begin{table}[tbp] \centering
%EndExpansion%
\begin{tabular}
[c]{|c|c|c|}\hline
$n$ & $\tau_{\ast}^{2}=1$ & $\tau_{\ast}^{2}=0.5$\\\hline
$5$ & $0.107$ & $0.193$\\
$10$ & $0.075$ & $0.146$\\
$20$ & $0.051$ & $0.107$\\
$50$ & $0.031$ & $0.067$\\
$100$ & $0.021$ & $0.046$\\\hline
\end{tabular}
\caption{Average bias against $H_0={0}$  when using a $N(0,\tau_*^2)$ prior for different sample sizes
$n$.}\label{loctab3}%
%TCIMACRO{\TeXButton{E}{\end{table}}}%
%BeginExpansion
\end{table}%
%EndExpansion

To use (\ref{biasagest2}) it is necessary to maximize $M(RB(\mu\,|\,X)\leq
1\,|\,\mu)$ as a function of $\mu$ and it is seen that, at least when the
prior is not overly concentrated, that this maximum occurs at $\mu=\mu_{\ast
}.$ When using the $N(0,1)$ prior the maximum occurs at $\mu=0$ when $n=5$ and
from the second column of Table \ref{loctab1}, the maximum equals $0.143$. The
average bias against is $0.107,$ as recorded in Table \ref{loctab3}.
\smallskip

Bias in favor occurs when the prior probability that $Im_{\Psi}$ does not
cover a false value is large, as this would seem to imply, ignoring the
situation where there is no evidence either way, that the plausible region
will cover a randomly selected false value from the prior with high prior
probability. Considering only those values that differ meaningfully from the
true value, the bias in favor for estimation is given by%
\begin{align}
&  E_{\Pi_{\Psi}}\left(  \sup_{\psi:d_{\Psi}(\psi,\psi_{0})\geq\delta}%
M(\psi_{0}\notin Im_{\Psi}(X)\,|\,\psi)\right)  \nonumber\\
&  =E_{\Pi_{\Psi}}\left(  \sup_{\psi:d_{\Psi}(\psi,\psi_{0})\geq\delta
}M(RB_{\Psi}(\psi_{0}\,|\,X)\geq1\,|\,\psi)\right)  ,\label{biasfavest2}%
\end{align}
which is the prior mean biases in favor for $H_{0}:\Psi(\theta)=\psi_{0}$
averaged over $\psi_{0}.$

\begin{figure}
[t]
\begin{center}
\includegraphics[
height=2.6593in,
width=2.6593in
]%
{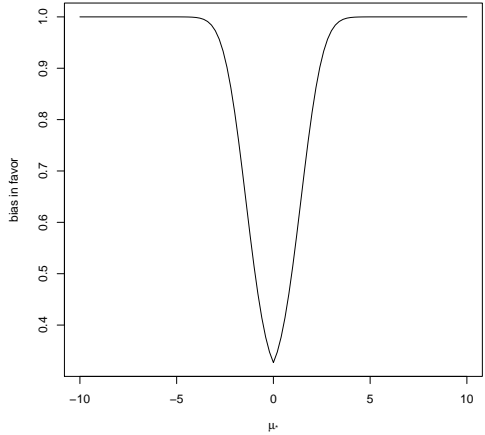}%
\caption{Bias in favor of\ $\mu$\ maximized over $\mu\pm\delta$ based on a
$N(0,1)$ prior with $\sigma_{0}=1,n=20,\delta=0.5.$}%
\label{fig3}%
\end{center}
\end{figure}

An upper bound on (\ref{biasfavest2}) is commonly equal to 1 as illustrated in
Figure \ref{fig3} and so is not useful.$\smallskip$

\noindent\textbf{Example 5. }\textit{Location normal (continued).}

Figure \ref{fig3} plots $\sup_{\,\mu\pm\delta}M(\mu\notin Im_{\Psi}%
(X)\,|\,\mu)$ as a function of $\mu$ for a specific prior and $\delta.$ Table
\ref{loctab5} contains values of (\ref{biasfavest2}) for this situation with
different values of $\delta$.%
%TCIMACRO{\FRAME{ftFU}{2.6593in}{2.6593in}{0pt}{\Qcb{Bias in favor of\ $\mu
%$\ maximized over $\mu\pm\delta$ based on a $N(0,1)$ prior with $\sigma
%_{0}=1,n=20,\delta=0.5.$}}{\Qlb{fig3}}{fig3.eps}%
%{\special{ language "Scientific Word";  type "GRAPHIC";  display "USEDEF";
%valid_file "F";  width 2.6593in;  height 2.6593in;  depth 0pt;
%original-width 6.9998in;  original-height 6.9998in;  cropleft "0";
%filename 'fig3.eps';file-properties "XNPEU";}}}%
%BeginExpansion
%EndExpansion%
%TCIMACRO{\TeXButton{B}{\begin{table}[tbp] \centering}}%
%BeginExpansion
\begin{table}[tbp] \centering
%EndExpansion%
\begin{tabular}
[c]{|c|c|c|}\hline
$n$ & $(\mu_{\ast},\tau_{\ast})=(0,1),\delta=1.0$ & $(\mu_{\ast},\tau_{\ast
})=(0,1),\delta=0.5$\\\hline
$5$ & $0.451$ & $0.798$\\
$10$ & $0.185$ & $0.690$\\
$20$ & $0.025$ & $0.486$\\
$50$ & $0.000$ & $0.131$\\
$100$ & $0.000$ & $0.009$\\\hline
\end{tabular}
\caption{Average bias in favor for estimation  when using a $N(0,\tau_*^2)$ prior for different sample sizes
$n$ and difference $\delta$. }\label{loctab5}%
%TCIMACRO{\TeXButton{E}{\end{table}}}%
%BeginExpansion
\end{table}%
%EndExpansion

\section{Properties}

The inferences possess a variety of good properties and are even are optimal
in the family of all Bayesian inferences. For example, the invariance of the
inference under reparameterizations has already been noted. Another simple
property is the relationship
\begin{equation}
RB_{\Psi}(\psi\,|\,x)=E_{\Pi(\cdot\,|\,\psi)}(RB_{\Psi}(\theta\,|\,x))
\label{av}%
\end{equation}
where $\Pi(\cdot\,|\,\psi)$ is the conditional prior for $\theta$ given
$\Psi(\theta)=\psi.$ So the evidence concerning $\psi$ is obtained by
averaging the evidence for each value of $\theta\in\Psi^{-1}\{\psi\}$ with
respect to the conditional distribution on the nuisance parameters. This
relationship leads to a result that may appear paradoxical as discussed in the
following example.$\smallskip$\newpage

\noindent\textbf{Example 6. }\textit{Are evidence measures incoherent?}

Consider $A,B,C$ subsets of $\Omega$ with probability measure $P$ such that
$P(A)>0,P(B)>0,P(C)>0$ and suppose $A\cap B=\phi.$ Then applying (\ref{av})
gives%
\[
RB(A\cup B\,|\,C)=RB(A\,|\,C)P(A\,|\,A\cup B)+RB(B\,|\,C)P(B\,|\,A\cup B).
\]
From this it can be deduced that when $RB(A\,|\,C)>1,$ then $RB(A\cup
B\,|\,C)<1$ iff $RB(B\,|\,C)<1$ and%
\begin{equation}
P(A\,|\,A\cup B)<\frac{1-RB(B\,|\,C)}{RB(A\,|\,C)-RB(B\,|\,C)}. \label{av2}%
\end{equation}
So it is possible that there is evidence in favor of $A$ but evidence against
the superset $A\cup B.$ The Bayes factor also has this property and this has
been referred to as \textquotedblleft incoherence\textquotedblright\ of a
measure of evidence.\ It is clear from (\ref{av2}) that $A$ has to be a
suitably small proportion of $A\cup B$ for this to hold. But note that, if
$C=A,$ then $RB(B\,|\,C)=0$ and, even if $P(A\,|\,A\cup B)$ is very small,
then (\ref{av2}) fails. So the situation is more subtle.

To see how such a situation can arise consider the following. In a population
$\Omega$ let $C$ correspond to a trait that exists in subpopulation $A\cup B$
only in $A.$ So $RB(A\,|\,C)=1/P(C)>1$ and $RB(B\,|\,C)=0.$ This implies
$RB(A\cup B\,|\,C)=P(A\,|\,A\cup B)/P(C)$ and so there is evidence against
$A\cup B$ iff $P(A\,|\,A\cup B)<P(C)$ or iff the probability of possessing the
trait within the subpopulation is smaller than the probability of possessing
the trait within the full population. This makes sense and so our claim is
that there is no incoherency in the relative belief ratio (or Bayes factor) as
a measure of evidence. Measuring evidence is seen once again to be different
than measuring belief.$\smallskip$

In Evans and Guo (2019) a number of results are established concerning the
principle of evidence. For example, the following result establishes the
unbiasedness of the inferences that result.

\noindent\textbf{Theorem 1.} Using the principle of evidence (i) the prior
probability of getting evidence in favor of $\psi_{0}$ given that it is true
is greater than or equal to the prior probability of getting evidence in favor
of $\psi_{0}$ given that $\psi_{0}$ is false and (ii) the prior probability of
$Pl_{\Psi}$\ covering the true value is always greater than or equal to the
prior probability of $Pl_{\Psi}$ covering a false value.\smallskip

\noindent Moreover the principle of evidence is seen to satisfy a number of
optimality properties. For example, suppose another principle is chosen for
establishing evidence against or in favor of a value $\psi$ and let
$D(\psi)\subset\mathcal{X}$ be the data sets that don't lead to evidence in
favor of $\psi.$ For the principle of evidence this is the set $R(\psi
)=\{x:RB_{\Psi}(\psi\,|\,x)\leq1\}.$ The rules that are potentially of
interest satisfy
\begin{equation}
M(D(\psi)\,|\,\psi)\leq M(R(\psi)\,|\,\psi) \label{ineq1}%
\end{equation}
as this indicates that $D(\psi)$ will do no worse than $R(\psi)$ at
identifying when $\psi$ is true and is better when the inequality is strict.
The following Theorem holds where part (ii) is the relevant result and part
(i) is used to obtain this.\smallskip

\noindent\textbf{Theorem 2.} (i) The prior probability $M(D(\psi))$ is
maximized among all rules satisfying (\ref{ineq1}) by $D(\psi)=R(\psi).$ (ii)
If $\Pi_{\Psi}(\{\psi\})=0,$ then $R(\psi)$ maximizes, among all rules
satisfying (\ref{ineq1}), the prior probability of not obtaining evidence in
favor of $\psi$ when it is false and otherwise maximizes this probability
among all rules satisfying $M(D(\psi)\,|\,\psi)=M(R(\psi)\,|\,\psi
).$\smallskip

\noindent Now consider $C(x)=\{\psi:x\notin D(\psi)\}$ which is the set of
$\psi$ values for which there is evidence in their favor after observing $x$
according to some alternative evidence rule. If (\ref{ineq1}) holds for each
$\psi,$ then $E_{\Pi_{\Psi}}\left(  M(\psi\in C(X)\,)\,|\,\psi)\right)  \geq
E_{\Pi_{\Psi}}\left(  M(\psi\in Pl_{\Psi}(X)\,)\,|\,\psi)\right)  $ and so the
Bayesian coverage of $C$ is at least as large as that of $Pl_{\Psi}$ and so
represents a viable alternative to using $Pl_{\Psi}.$ The following
establishes an optimality result for $Pl_{\Psi}$ and again it is part (ii)
that is relevant.\smallskip

\noindent\textbf{Theorem 3.} (i) The prior probability that the region $C$
doesn't cover a value $\psi$ generated from the prior, namely, $E_{\Pi_{\Psi}%
}(M(\psi\notin C(X))),$ is maximized, among all rules satisfying (\ref{ineq1})
for all $\psi,$ by $C=Pl_{\Psi}.$ (ii) If $\Pi_{\Psi}(\{\psi\})=0$ for all
$\psi,$ then $Pl_{\Psi}$ maximizes, among all rules satisfying (\ref{ineq1})
for all $\psi,$ the prior probability of not covering a false value and
otherwise maximizes this probability among all $C$ satisfying $M(\psi\notin
C(X)\,|\,\psi)=M(\psi\notin Pl_{\Psi}(X)\,|\,\psi)$ for all $\psi.$\smallskip

\noindent If there is a value $\psi_{\ast}=\arg\inf_{\psi}M(\psi\in Pl_{\Psi
}(X)\,)\,|\,\psi),$ then $\gamma_{\ast}=M(\psi_{\symbol{94}}\in Pl_{\Psi
}(X)\,)\,|\,\psi_{0})$ serves as a lower bound on the coverage probabilities,
and thus $Pl_{\Psi}$ is a $\gamma_{\ast}$-confidence region for $\psi$ and
this is a pure frequentist $\gamma_{\ast}$-confidence region when $\Psi
(\theta)=\theta.$ Since $M(\psi\in Pl_{\Psi}(X)\,)\,|\,\psi)=1-M(\psi\notin
Pl_{\Psi}(X)\,)\,|\,\psi)=1-M(R(\psi)\,|\,\psi),$ then Example 1 shows that it
is reasonable to expect that such a $\psi_{0}$ exists. Optimality results
similar to Theorems 2 and 3 also exist when considering evidence in favor and
are presented in Evans and Guo (2019).

As discussed in Evans (2015) relative belief credible regions $C_{\Psi,\gamma
}$ are unbiased and possess optimality properties similar to those stated for
$Pl_{\Psi}$.\ Also, it is always the case that $RB(C_{\Psi,\gamma
}(x)\,|\,x)\geq1,$ with the inequality generally strict. This implies that
there is evidence that the true value is in $C_{\Psi,\gamma}(x)$ and this is
not the case with other rules for forming credible regions. Also, it can be
proved that $RB(C\,|\,x)$ is maximized by $C=C_{\Psi,\gamma}(x)$ among all
sets $C$ satisfying $\Pi_{\Psi}(C\,|\,x)=\Pi_{\Psi}(C_{\Psi,\gamma
}(x)\,|\,x).$ So $C_{\Psi,\gamma}$ represents those values that have the
greatest increase in belief amongst all the possible $\gamma$ credible regions
having the same posterior content.

Evans (2015) contains discussion concerning a number of other properties for
these inferences. Overall the inferences derived via R$_{\text{2}}$ and
R$_{\text{3}}$ are seen to have excellent properties.

\section{Conclusions}

This paper has described an approach to statistical reasoning. The central
idea, which in fact we view as key for any such theory, is to be precise about
the characterization of statistical evidence relevant to a particular problem.
This is accomplished via the principle of evidence which overall plays the
major role. This principle is intuitively very appealing and, while very
simple, has broad implications. The principle has been extensively discussed
in the philosophy of science literature but statisticians have for the most
part ignored it, and that may well be because of the requirement that prior
beliefs be specified via a probability measure.

Another feature of the approach is a clear separation between the
specification and checking of the ingredients and the rules of inference. Two
necessary ingredients are a statistical model and a prior. Both of these are
based upon choices made by an analyst, perhaps in conjunction with application
area experts. As such, these ingredients are subjective in nature and it is
fair to question their validity in light of the central role that the goal of
objectivity plays in science. This is addressed in several ways. One is to
acknowledge that the primary role of these ingredients is to allow for the
reasoning process to present answers to the questions of scientific interest
about the real-world object $\Psi.$ In essence the correctness of the
ingredients is not of primary importance and obsessing about this seems
misplaced. Still, it is clear that the ingredients can be chosen very badly in
the sense that they are strongly contradicted by the data and this may lead to
misleading inferences. But there is a scientifically valid approach to dealing
with this issue through checking the model and prior against the objective, at
least if it is collected properly, data. This may lead to the need to modify
chosen ingredients and, at least for the prior, this can be done in a way that
is only very weakly data dependent.

The chosen ingredients can also suffer from an even more serious defect as the
choices may seriously bias the inferences. This bias is measured by the prior
probabilities associated with obtaining evidence in favor of or against. The
principle of evidence is seen to play the central role in measuring and
controlling these biases and can even be considered the optimal way to do
this. For us, the measurement of bias is an absolute necessity in assessing
the value of a statistical study. If there is a high prior probability of
obtaining misleading evidence in a study, then the value of any evidence that
reflects this is questionable. Another consequence of measuring and
controlling bias is that it establishes a cooperative and beneficial
relationship between frequentist and Bayesian ideas, as one can't do without
the other. Frequentism is seen to be an a priori concern with controlling the
inferences that can result from potential data sets, while Bayesianism is
concerned with the actual inferences that follow from the rules of inference
applied to the specified ingredients together with the observed data.

One counter-argument to what has been presented here is that it is very
idealistic, at least compared with many of the problems statisticians commonly
face. Our response to this is that the theory is definitely focused on
problems where, for example, the data collection process can be designed so
that appropriate and adequate amounts of data can be collected, where
parameters indexing models are sufficiently well-understood so that priors can
be elicited, differences of importance (values like $\delta$) can be specified
and are such that the computations can be carried out to sufficient accuracy.
Such problems are what might be considered the core statistical problems. If
an appropriate theory of statistical reasoning cannot be developed for such
contexts, then there is little reason to hope for a theory that can handle the
kinds of problems where data collection is not designed, the data are far too
few for the dimensionality of the problem and aspects of the models used are
not well understood so that default choices of priors are employed. In any
case a proposed theory has to handle the core problems and deal with the
concept of statistical evidence appropriately. A better approach, in our view,
is to identify the problems which can be considered central and, as discussed
in the paper, when problems are encountered that lie outside this core
consider compromises to the theory that retain as much of the ideal approach
as possible. Certainly Cao et al. (2015), and Evans and Tomal (2018) are in
this spirit. Additional applications to more significant statistical problems
than those discussed here can be found in the papers cited in this paper. Of
some significance are recent applications in quantum state estimation that
arises out of an approach to inference developed by B. Englert, see Shang et
al. (2013), which has much in common with what has been presented here. Gu et
al. (2019) is an example of this and is one where the concept of bias played a
notable role.

The computations associated with the theory can indeed be difficult. For the
computations associated with bias and prior-data conflict, however, this is
mitigated by requiring only a very low level of accuracy for the relevant
probabilities in question. See, for example, Nott et al. (2016), Nott et al.
(2019) and Wang et al. (2018) where approximate calculation approaches have
been successfully employed. In any case, we subscribe to the view that is
better to approximately compute what is believed to be correct rather than
exactly compute what isn't. It is notable that most of the hard work involved
in an application is associated with choosing and checking the ingredients
and, since this is the aspect that involves the most thinking about the
application and the relevance of the statistical analysis to it, this is as it
should be.

\section{References}

\noindent Achinstein, P. (2001) The Book of Evidence. Oxford University
Press.\smallskip

\noindent Aitkin, M. (2010) Statistical Inference: An Integrated
Bayesian/Likelihood Approach. Chapman and Hall/CRC.\smallskip

\noindent Al-Labadi, L. and Evans, M. (2017) Optimal robustness results for
some Bayesian procedures and the relationship to prior-data conflict. Bayesian
Analysis 12, 3, 702-728.\smallskip

\noindent Al-Labadi, L., Baskurt, Z and Evans, M. (2017) Goodness of fit for
the logistic regression model using relative belief. J. of Statistical
Distributions and Applications, 4:17.\smallskip

\noindent Baskurt, Z. and Evans, M. (2013) Hypothesis assessment and
inequalities for Bayes factors and relative belief ratios. Bayesian Analysis,
8, 3, 569-590.\smallskip

\noindent Berger, J.O. and Selke, T. (1987) Testing a point null hypothesis:
the irreconcilability of p values and evidence. Journal of the American
Statistical Association, 82, 397, 112-122.\smallskip

\noindent Berger, J.O. and Delampady, M. (1987) Testing precise hypotheses.
Statistical Science, 2, 3, 317-335.\smallskip

\noindent Birnbaum, A. (1964) The Anomalous Concept of Statistical Evidence:
Axioms, Interpretations, and Elementary Exposition. New York University,
Courant Institute of Mathematical\ Sciences.\smallskip

\noindent Box, G. (1980). Sampling and Bayes' inference in scientific
modelling and robustness. Journal of the Royal Statistical Society, A, 143:
383--430.\smallskip

\noindent Brown, L. D., Cai, T., DasGupta, A. (2002) Confidence Intervals for
a binomial proportion and asymptotic expansions. Ann. Statist. 30, 1,
160--201.\smallskip

\noindent Cao, Y., Evans, M. and Guttman, I. (2015) Bayesian factor analysis
via concentration. Current Trends in Bayesian Methodology with Applications,
edited by S. K. Upadhyay, U. Singh, D. K. Dey and A. Loganathan, 181-201, CRC
Press.\smallskip

\noindent Evans, M. (2015) Measuring Statistical Evidence Using Relative
Belief. Chapman and Hall/CRC.\smallskip

\noindent Evans, M. and Guo, Y. (2019) Measuring and controlling bias for some
Bayesian inferences and the relation to frequentist criteria.
arXiv:1903.01696\smallskip

\noindent Evans, M., Guttman, I. and Li, P. (2017) Prior elicitation,
assessment and inference with a Dirichlet prior. Entropy 2017, 19(10),
564.\smallskip

\noindent Evans, M. and Jang, G-H. (2011a) A limit result for the prior
predictive applied to checking for prior-data conflict. Statistics and
Probability Letters, 81, 1034-1038.\smallskip

\noindent Evans, M. and Jang, G-H. (2011b) Weak informativity and the
information in one prior relative to another. Statistical Science, 26, 3,
423-439.\smallskip

\noindent Evans, M. and Moshonov, H. (2006) Checking for prior-data conflict.
Bayesian Analysis, 1, 4, 893-914.\smallskip

\noindent Evans, M. and Tomal, J. (2018) Multiple testing via relative belief
ratios. FACETS, 3: 563-583, DOI: 10.1139/facets-2017-0121.\smallskip

\noindent Kneale, W. and Kneale, M. (1962) The Development of Logic. Clarendon
Press.\smallskip

\noindent Gelman, A., Jahukin, A., Pittau, M. G. and Su, Y.-S. (2008) A weakly
informative default prior distribution for logistic and other regression
models. Ann. Appl. Statist. 2 1360--1383.\smallskip

\noindent Gu, Y., Li, W. Evans, M. and Englert, B-G. (2019) Very strong
evidence in favor of quantum mechanics and against local hidden variables from
a Bayesian analysis. Physical Review A 99, 022112 (1-17).\smallskip

\noindent Morey, R., Romeijn, J-W, and Rouder, J. (2016) The philosophy of
Bayes factors and the quantification of statistical evidence. J. of
Mathematical Psychology, 72, 6-18.\smallskip

\noindent Mosteller, F. and Tukey, J. (1977). Data analysis and regression : a
second course in statistics. Addison-Wesley.\smallskip

\noindent Nott, D., Drovandi, C., Mengersen, K. and Evans, M. (2016)
Approximation of Bayesian predictive p-values with regression ABC. Bayesian
Analysis 13, 1, 59-83.\smallskip

\noindent Nott, D.J., Seah, M., Al-Labadi, L., Evans, M., Ng, H.K., and
Englert, B-G. (2019) Using prior expansions for prior-data conflict checking.
arXiv:1902.10393.\smallskip

\noindent Nott,D., Wang, X., Evans, M., and Englert, B-G. (2018) Checking for
prior-data conflict using prior to posterior divergences.
arXiv:1611.00113.\smallskip

\noindent Ripamonti, E., Lloyd, C., and Quatto, P. (2017) Contemporary
frequentist views of the 2x2 binomial trial. Statist. Sci. 32, 4,
600--615.\smallskip

\noindent Robert, C. P. (2014) On the Jeffreys-Lindley paradox. Philosophy of
Science, 81, 216--232.\smallskip

\noindent Royall, R. (1997) Statistical Evidence: A Likelihood Paradigm.
Chapman and Hall/CRC.\smallskip

\noindent Salmon, W. (1973) Confirmation. Scientific American, 228, 5,
75-81.\smallskip

\noindent Shafer, G. (1976) A Mathematical Theory of Evidence. Princeton
University Press.\smallskip

\noindent Shang, J., Ng, H. K., Sehrawat, A. Li, X. and Englert, B.-G. (2013)
Optimal error regions for quantum state estimation. New J. Phys. 15,
123026.\smallskip

\noindent Thompson, B. (2007) The Nature of Statistical Evidence. Lecture
Notes in Statistics 189, Springer.\smallskip

\noindent Vieland, V.J and Seok, S-J. (2016)\ Statistical evidence measured on
a properly calibrated scale for multinomial hypothesis comparisons. Entropy
18(4): 114.\smallskip

\noindent Villa, C. and Walker, S. (2017) On the mathematics of the
Jeffreys--Lindley paradox. Communications in Statistics - Theory and Methods,
46, 24, 12290-12298.\smallskip

\noindent Wang, X., Nott, D. J., Drovandi, C. C., Mengersen, K. and Evans, M.
(2018) Using History Matching for Prior Choice, Technometrics, 60:4, 445-460.

\end{document}